\documentclass[10pt,fleqn]{article}

\usepackage{amsfonts,amssymb,amsmath}
\usepackage{graphicx} 
\usepackage{color}
\usepackage[bookmarksnumbered, plainpages, backref, colorlinks]{hyperref}

\newcommand{\R}{\mathbb R}
\newcommand{\N}{\mathbb N}

\newcommand{\ds}{\displaystyle}
\newcommand{\lan}{\langle}
\newcommand{\ran}{\rangle}

\newcommand{\ul}{\underline}
\newcommand{\ol}{\overline}

\newcommand{\rha}{\rightharpoonup}

\newcommand{\qed}{\hfill $\Box$}

 \def\reff#1{\mbox{\rm(\ref{#1})}}

\newtheorem{theorem}{Theorem}
\newtheorem{lemma}{Lemma}
\newtheorem{proposition}{Proposition}

\def\div{\mbox{\rm div\hspace{1mm}}}

\def\hnorm#1#2{\vert\,#1\,\vert_{#2}}

\title{On optimal control in a  nonlinear interface problem  
 described by hemivariational inequalities }
\author{Joachim Gwinner  \\ 
 Department of Aerospace Engineering \\
 Universit\"{a}t der Bundeswehr M\"{u}nchen \\ 
Neubiberg - Munich, Germany} 

\begin{document}
\date {} 
\maketitle

\begin{center} 
{\it Dedicated to Professor Thomas Apel \\
on the occasion of his 60th birthday}
\end{center}

\begin{abstract}
The purpose of this paper is three-fold.
Firstly we attack a nonlinear interface problem on an unbounded domain  
with nonmonotone set-valued transmission conditions.
The investigated problem involves a nonlinear monotone partial differential 
equation in the interior domain and the Laplacian in the exterior domain.
Such a scalar interface problem models nonmonotone frictional contact of 
elastic infinite media. The variational formulation of the interface problem
leads to a hemivariational inequality (HVI), which however lives on the
unbounded domain, and thus cannot analyzed in a reflexive Banach space setting.
By boundary integral methods we  obtain another HVI that is amenable to 
functional analytic methods using standard Sobolev spaces on the interior
domain and Sobolev spaces of fractional order on the coupling boundary.
Secondly broadening the scope of the paper, we consider extended 
real-valued HVIs augmented by convex extended real-valued functions. Under a
smallness hypothesis, we provide existence and uniqueness results, also
establish  a stability result with respect to the extended real-valued function as parameter.
Thirdly based on the latter stability result, we prove the existence of optimal controls for four kinds of optimal control problems: distributed control on the bounded domain, boundary control,  
simultaneous distributed-boundary control
 governed by the interface problem, as well as control of the obstacle driven by a related bilateral obstacle interface problem. 
\end{abstract}

\\[1ex]

\section{Introduction}

Optimal control of partial differential equations (PDEs) is a vast field of applied mathematics. Here we focus to the control of elliptic PDEs,
governed by a hemivariational inequality (HVI) in a weak formulation.   
 
The theory of HVIs was introduced and has been 
studied since 1980s by Panagiotopoulos \cite{Pan-1993}, as a generalization of
variational inequalities with the aim to model many problems coming from
mechanics when the energy functionals are nonconvex, but locally Lipschitz,
so the Clarke generalized differentiation calculus \cite{Clarke} 
can be used, see \cite{Naniewicz,GMDR-2003,GoeMot-2003}.
For more recent monographs on HVIs with application to
contact problems we refer to \cite{Ochal,SofMig-2018}. 

While optimal control in variational inequalities has already been 
treated for a longer time in various works,
see the monograph \cite{Ba1984} 
and e.g. the papers \cite{Pa1977,MiPu1984,Fr1986,ALY1998,ItKu2000,DLR2011},
optimal control in HVIs has been more recently studied, see e.g.
\cite{HaPa1995,PeKu2018,Sof2018,LiMi2019,Sof2019,SoMa2019,LYZe2021}.
Contrary to the cited work above the underlying state problem of this paper is not a boundary value problem on a bounded domain, but an interface problem involving a pde on an unbounded domain. For the simplicity of presentation  we discuss a scalar interface problem with a monotone pde on the interior domain, the Laplacian on the exterior domain, and - a main novelty of the paper -
nonmonotone set-valued transmission conditions.
This scalar problem models nonlinear contact problems with nonmonotone friction in infinite elastic media that arise in various fields of science and technology; let us mention geophysics, see e.g. \cite{scholz2019}, soil mechanics, in particular soil-structure interaction problems, see e.g. \cite{GhIo2009}, and civil engineering of underground structures, see e.g. \cite{wang2018soil}.

It should be underlined that such interface problems involving a pde on an unbounded domain are more difficult than standard boundary value problems on bounded domains, since  a direct variational formulation 
of the former problems leads to a HVI, which lives on the
unbounded domain, and thus cannot analyzed in a reflexive Banach space setting.
At this point boundary integral methods, see the monograph \cite{HsWe-2008},
come into play. By these methods we  obtain another HVI that is amenable to 
functional analytic methods using standard Sobolev spaces on the interior
domain and Sobolev spaces of fractional order on the coupling boundary. 
Let us note in passing that these integral methods lay the basis for the numerical treatment of such interface problems by the well-known coupling
of boundary elements and finite elements, 
see \cite[Chapter 12]{GwiSte-2018}.

A further novel ingredient of our analysis is a stability theorem
that considerably improves a related result in the recent paper 
\cite{Sof2018} and extends it to more general extended 
real-valued HVIs augmented by convex extended real-valued functions.
This stability theorem provides the key to a unified approach to 
the existence of optimal controls in various optimal control problems:
 distributed control on the bounded domain, boundary control,
simultaneous distributed-boundary control
  governed by the interface problem,
 as well as the control of the obstacle in a related bilateral obstacle
interface problem.    

The plan of the paper is as follows. 
The next section 2 consists of three parts:
a collection of some basic tools of Clarke's generalized
differential calculus for the analysis of the nonmonotone
transmission conditions, 
a description of the interface problem under study in strong form, and
existence and uniqueness results for a class of abstract HVIs using an equilibrium approach.  
In section 3 a first equivalent weak variational formulation of the
interface problem is derived in terms of a HVI, which lives on the unbounded domain (as the original problem). Section 4 employs boundary integral
methods and provides a HVI formulation of the interface problem on the
bounded boundary-domain. Section 5 establishes well-posedness results,
in particular a stability theorem for a more general class of extended
real-valued HVIs. Based on this stability theorem, section 6 presents 
a unified approach to the existence of optimal controls in four  optimal control problems: distributed, boundary, boundary-distributed and obstacle control. 
The final section 7 shortly summarizes our findings, gives some concluding
remarks, and sketches some directions of further research.

\section {Some preliminaries, the interface problem, and an equilibrium
approach to HVIs} \label{inter}
\setcounter{equation}{0}
\subsection{Some preliminaries from Clarke's generalized differential calculus}

Before we pose our interface problem, let us first recall the central notions
of Clarke's generalized differential calculus \cite{Clarke}.
Let $\phi: X \to \R$ be a locally Lipschitz function 
on a real Banach space $X$. Then  
$$ \phi^0(x;z) := \ds \limsup_{ y \to x; t \downarrow 0}
  \frac{\phi(y+tz) - \phi(y)}{t} ~ ~ x, z  \in X,  $$
is called the {\it generalized directional derivative} of $\phi$ 
at $x$ in the direction $z$.
Note that the function $z \in X \mapsto \phi^0(x;z)$ is finite, positively
homogeneous, and sublinear, hence convex and continuous; further, the function
$(x,z) \mapsto \phi^0(x;z)$ is upper semicontinuous.
The {\it generalized gradient} of the function $\phi$ at $x$,
 denoted by (simply) $\partial \phi(x)$,
is the unique nonempty $\text{weak}^*$ compact convex subset of
the dual space $X'$, whose support function is $\phi^0(x;.)$.
Thus 
\begin{eqnarray*} 
&& \xi \in \partial \phi(x) \Leftrightarrow \phi^0(x;z) 
\ge \lan \xi,z \ran, \, \forall
z \in X , \\[0.5ex]
&& \phi^0(x;z) = \max \{ \lan \xi,z \ran ~:~ \xi \in \partial \phi(x) \}, \, \forall z
\in X \,.
\end{eqnarray*}
When $X$ is finite dimensional, then, according to Rademacher's theorem,
 $\phi$ is differentiable almost everywhere, and the generalized gradient of $\phi$ at a point $x\in \R^n$ can be characterized by
$$
\partial \phi(x)= \textnormal{co} \,
 \{\xi\in \R^{ n} \, : \, \xi = \displaystyle \lim _{k \to \infty}
\nabla \phi(x_k), \; x_k \to x, \,\phi \, \textnormal{ is differentiable at }\, x_k\},$$
where "co" denotes  the convex hull.
\subsection{The interface problem - strong formulation}
Let $\Omega \subset \R^d ~(d \ge 2)$ be a bounded domain with Lipschitz
 boundary $\Gamma$. 
% In the following we focus to the case $d =3$; 
% in the case $d=2$ some pecularities of the
% analysis by boundary integral methods come up that need extra attention, 
% see e.g.  \cite[chapter 12]{GwiSte-2018}.
 To describe mixed transmission conditions, we  
have  two non-empty, open  disjoint boundary parts 
$ \Gamma_s$ and $\Gamma_t$ such that the boundary $\Gamma = 
\textnormal{cl  } \Gamma_s \cup \textnormal{cl } \Gamma_t$.
Let $n$ denote the unit normal on $\Gamma$
defined almost everywhere
pointing from $\Omega $ into $\Omega^c :=\R^d\setminus\overline{\Omega} $.
Let the data $f\in L^2(\Omega)$, $u_0\in H^{1/2}(\Gamma)$, and
$q \in L^2(\Gamma)$ be given.

In the interior part $\Omega$, we consider the nonlinear partial differential
equation
\begin{equation}\label{a1}
\div \Bigl( p(| \nabla u |)\cdot \nabla u \Bigr) + f = 0 
\qquad \mbox{in } \; \Omega,
\end{equation}
where $p:[0,\infty) \to [0,\infty)$ is a continuous
function with $ t \cdot p(t) $
being  monotonously increasing with $t$.

In the exterior part  $\Omega^c$, we consider the Laplace equation
\begin{equation}\label{a2}
\Delta u=0 \qquad \mbox{in }\; \Omega^c
\end{equation}
with the radiation condition at infinity  for $|x|\to\infty$,
\begin{eqnarray}
\label{a3}     %\nonumber
u(x) =  
    \left\{ \begin{array}{ll}
     a + o(1)  & \mbox{if } d=2 \,, \\[0.5ex]  
      O(|x|^{2-d}  )  & \mbox{if } d > 2 \,, 
    \end{array} \right\}
  \end{eqnarray}
where $a$ is a real constant for any $u$, but may vary with $u$.

Let us write $u_1:=u|_\Omega $ and $u_2:=u|_{\Omega^c} $, then the tractions on
the coupling boundary $\Gamma$ are given by the traces of
$ p(|\nabla u_1|) \frac{\partial u_1}{\partial n}  $ and
$-\frac{\partial u_2}{\partial n} $, respectively. 
 
We prescribe classical transmission conditions on $\Gamma_t$,
\begin{equation}\label{a4}
u_1|_{\Gamma_t} = u_2|_{\Gamma_t}+ u_0|_{\Gamma_t} \quad\mbox{and}\quad
p(|\nabla u_1|)\left.\frac{\partial u_1}{\partial n}\right \vert_{\Gamma_t}
= \left.\frac{\partial u_2}{\partial n}\right\vert_{\Gamma_t} + q|_{\Gamma_t}
,
\end{equation}
and  on $\Gamma_s$, analogously for the tractions, 
\begin{equation} \label{a5_1}
p(|\nabla u_1|)\left.\frac{\partial u_1}{\partial n}\right\vert_{\Gamma_s}
= \left.\frac{\partial u_2}{\partial n}\right\vert_{\Gamma_s} + q|_{\Gamma_s}
\end{equation}
and the generally nonmonotone, set-valued transmission condition,   
\begin{equation}  \label{a5}
p(|\nabla u_1|) \left. \frac{\partial u_1}{\partial n}\right\vert_{\Gamma_s}  \in \left. \partial j(\cdot, u_0  + (u_2 - u_1)\right\vert_{\Gamma_s}) 
 \,.
\end{equation}

Here the function $j: \Gamma_s \times \R \to \R $ is such that $j(\cdot,
\xi):\Gamma_s\to \R$ is measurable on $\Gamma_s$ for all $\xi \in \R$ and $j(s,
\cdot) : \R \to \R$ is locally Lipschitz for almost all (a.a.) $s\in \Gamma_s$
with
$\partial j(s,\xi):=\partial j(s,\cdot) (\xi)$, the generalized gradient of
$j(s,\cdot)$ at $\xi$.

Further, we require the following growth condition on the 
so-called superpotential $j$:
There exist positive constants $c_{j,1}$ and $c_{j,2}$ such that for a.a. $s \in
\Gamma_s$, all $\xi \in \R$ and for all
$\eta \in \partial j(s, \xi)$ the following inequalities hold
\begin{equation} \label{as-j} % (H(j)) 
          (i)\quad |\eta| \leq c_{j,1} (1+|\xi|) \,, 
          (ii)\quad \eta ~\xi \geq - c_{j,2} |\xi| \,.
     \end{equation}
 
 Thus  the interface problem under study consists in finding   
 $u_1\in H^1(\Omega)$ and $u_2\in H^1_{loc}(\Omega^c )$
that satisfy  \reff{a1}--\reff{a5} in a weak form.

\subsection{An equilibrium approach to a class of HVIs - 
existence and uniqueness results}

Next we  describe the  functional analytic setting
for the interface problem and provide existence and uniqueness results
using an equilibrium approach.
To this end, let $X:= L^2 (\Gamma_s)$ and introduce the real-valued locally Lipschitz functional
\begin{equation} \label{def-J}
J(y) := \int_{\Gamma_s} j(s,y(s)) ~ ds\,,  \qquad y \in X \,.
\end{equation}
Then by Lebesgue's theorem of majorized convergence,
\begin{equation} \label{def-J0}
J^0(y;z)= \int_{\Gamma_s} j^0(s,y(s);z(s)) ~ ds\,, 
 \qquad (y,z) \in X\times X  \,,
\end{equation}
where $j^0(s, \cdot~;~\cdot)$ denotes the generalized directional derivative of
$j(s,\cdot)$.

As we shall see in the subsequent sections, the weak formulation of the problem
\reff{a1}--\reff{a5} leads, in an abstract setting, to a hemivariational
inequality (HVI) with a nonlinear operator ${\cal A}$ and the nonsmooth
functional $J$, namely, we are looking for some
$\hat v\in {\cal C} $ such that
\begin{equation}\label{a6}
{\cal A}(\hat v) (v- \hat v) 
+ J^0(\gamma \hat v; \gamma v- \gamma \hat v)  %\varphi(\hat{v},v)
\ge \lambda(v-\hat v)
\qquad \forall v \in {\cal C}.
\end{equation}
Here ${\cal C}\ne \emptyset $ is a closed convex subset of a real reflexive Banach space $E$, the nonlinear operator ${\cal A} : E\to E^* $ is a monotone operator,  $\gamma := \gamma_{E \to X}$ is a linear continuous operator,
and the linear form $\lambda$ belongs to the dual $E^*$.
Similar to \cite{CCJG-1997}, the
operator ${\cal A}$ consists of a nonlinear monotone differential operator (as
made precise below) that results from the PDE \reff{a1} in the bounded domain
$\Omega$ and the positive definite Poincar\'e--Steklov operator on the boundary
$\Gamma$ of $\Omega$ that stems from
the exterior problem \reff{a2}-\reff{a3} and can be represented by the boundary
integral operators of potential theory. Thus it results that the operator
${\cal A}$ is Lipschitz continuous and strongly monotone with some monotonicity
constant $c_{\cal A} > 0$. 
On the other hand, by the compactness of the operator $\gamma$, which is seen below, the real-valued upper semicontinuous 
bivariate function, shortly bifunction  
$$\psi(v,w) :=  J^0(\gamma v; \gamma w- \gamma v)
  \,, \forall (v,w) \in E\times E  $$
becomes pseudomonotone, see 
\cite[Lemma 1]{OvGw-Rassias-14}.%, \cite[Lemma 4.1]{OvGw-2014}.
The latter result also shows a linear growth of $\psi(\cdot,0)$.
This and the strong monotonicity of ${\cal A}$ imply coercivity. 
Therefore by the theory of pseudomonotone VIs
\cite[Theorem 3]{Gwi-81}, \cite{Zei-2}, see \cite{OvGw-2014} 
for the application to HVIs, we have solvability  of  \reff{a6}. \\[1ex]
Further suppose that the generalized directional derivative $J^0$
satisfies the  one-sided Lipschitz condition: There exists $c_J > 0$
such that 
\begin{equation}\label{cc1}
J^0(y_1; y_2 - y_1) + J^0(y_2; y_1 - y_2)
\le c_J \| y_1 - y_2 \|_X^2 
\qquad \forall y_1,y_2 \in X \,.
\end{equation}
Then the smallness condition
\begin{equation}\label{cc2}
c_J \| \gamma \|^2_{E \to X} < c_{\cal A}
\end{equation}
 implies unique solvability of  \reff{a6}, see e.g. 
\cite[Theorem 5.1]{Ov2017} and \cite[Theorem 83]{SofMig-2018}.

It is noteworthy that under the smallness condition  \reff{cc2}
together with  \reff{cc1}, fixed point arguments 
\cite{Capa-2014} or the theory of set-valued pseudomonotone operators
\cite{SofMig-2018} are not needed, but simpler monotonicity arguments are
sufficient to conclude unique solvability. Moreover the compactness of the 
linear operator $\gamma$ is not needed either.
In fact, \reff{a6} can be framed as a
{\it monotone equilibrium problem} in the sense of Blum-Oettli \cite{BlOe94}:

\begin{proposition} \label{prop1}
 Suppose \reff{cc1} and \reff{cc2}. Then the bifunction
 $\varphi: {\cal C} \times {\cal C} \to \R$ 
defined by
\begin{equation}\label{bif1} 
\varphi(v,w) := 
{\cal A}(v) (w- v) 
+ J^0(\gamma v; \gamma w - \gamma v) - \lambda(w-v) 
\end{equation}
has the following properties: 
\begin{trivlist}
\item
\quad $\varphi(v,v) = 0 \,$ for all $v \in {\cal C}$;
 \item
\quad $\varphi(v,\cdot)$ is convex and lower semicontinuous 
for all $v \in {\cal C}$;
\item
\quad there exists some $\mu > 0$ such that 
$\varphi(v,w) + \varphi(w,v) \le  
- \mu  \|v - w \|^2_E \,$ for all \\ 
 $v,w \in {\cal C}$ (strong monotonicity);
\item
\quad the function $t \in [0,1] \mapsto \varphi(tw + (1-t)v,w)$ is upper
semicontinous at $t= 0$ for all $v,w \in {\cal C}$ (hemicontinuity).
\end{trivlist}

\end{proposition}
{\bf Proof.}
Obviously $\varphi$ vanishes on the diagonal and is convex and lower
semicontinuous with respect to the second variable.
To show strong monotonicity, estimate
\begin{eqnarray*}
&& \varphi(v,w) + \varphi(w,v) \\
&=& ({\cal A}(v) - {\cal A}(w)) (w- v) \\ 
&& + \, J^0(\gamma v; \gamma w - \gamma v) + J^0(\gamma w; \gamma v - \gamma w)
\\
&\le & - c_{\cal A} \, \|v - w \|^2_E + 
c_J \| \gamma v - \gamma w \|_X^2  \\
&\le & - (c_{\cal A} - c_J  \| \gamma \|^2_{E \to X} ) \, \|v - w \|^2_E \,. 
\end{eqnarray*}
To show hemicontinuity, it is enough to consider the 
bifunction $ (y,z) \in X \times X \mapsto J^0(y; z - y)$. Then for 
 $ (y,z) \in X \times X $ fixed, $t \in [0,1]$ one has
$$
J^0(y + t(z-y); z - (y + t(z-y))) = (1-t) J^0(y + t(z-y);z-y)
$$
and thus hemicontinuity follows from upper semicontinuity of $J^0$,
$$
\limsup_{t \downarrow 0}  J^0(y + t(z-y);z-y) \le  J^0(y;z-y) \,.
$$    
\qed

Since strong monotonicity implies coercivity and uniqueness, the fundamental
existence result \cite[Theorem 1]{BlOe94} applies to the HVI (\ref{a6}) to
conclude the following

\begin{theorem} \label{theo-1}
 Suppose \reff{cc1} and \reff{cc2}.
Then the  HVI (\ref{a6}) is uniquely solvable.
\end{theorem}

 \section{An intermediate HVI % hemivariational
 formulation of the interface problem}
\setcounter{equation}{0}
In this section we provide a first equivalent weak variational formulation of
the interface problem \reff{a1}~-~\reff{a5}
in terms of a hemivariational inequality (HVI). Since this 
HVI lives on the unbounded domain $\Omega \times \Omega^c$ 
(as the original problem), this HVI cannot treated 
in a reflexice Banch space setting  and therefore provides only 
an intermediate step in the analysis.

For the bounded Lipschitz domain $\Omega$ we use the standard Sobolev space
$H^s(\Omega)$
and the Sobolev spaces on the bounded Lipschitz boundary $\Gamma$  
(see \cite[Sect 2.4.1]{SaSch2011}),
\begin{eqnarray*} H^s(\Gamma) = 
    \left\{ \begin{array}{ll}
        \{ u|_{\Gamma} : u\in  H^{s+1/2}(\R^d) \} & (0 < s  \le 1) ,  \\
         L^2(\Gamma) & (s = 0),  \\
        (H^{-s}(\Gamma))^* \text{ (dual space) } & (-1 \le s <0).
    \end{array} \right.
  \end{eqnarray*}
Further we need for the unbounded domain $\Omega^c=\R^d \backslash
\overline{\Omega}$ the Frechet space
 (see e.g. \cite[Section 4.1, (4.1.43)]{HsWe-2008})
$$ H^s_{loc}(\Omega^c)  = \{u \in {\cal D }^*(\Omega^c) :
 \chi u \in H^s(\R^d) \,\, \forall \chi \in C_0^\infty (\Omega^c) \} \,.
$$ 

By the trace theorem we have $u|_{\Gamma}\in H^{1/2}(\Gamma)$ 
for $u \in H^1_{loc}(\Omega^c)$.
Next we define
$\Phi : H^1(\Omega)\times H^1_{loc}(\Omega^c) \to \R \cup\{\infty\} $
by
\begin{equation}\label{b1}
\Phi(u_1,u_2) := \int_{\Omega} g (|\nabla u_1|)\, dx +
                 \frac 12 \int_{\Omega^c} |\nabla u_2|^2\, dx
                - L(u_1,u_2|_\Gamma ).
\end{equation}

Here the data $f \in L^2(\Omega), q \in L^2(\Gamma)$
enter the linear functional
\begin{equation}\label{b2}
L(u,v) := \int_{\Omega} f \cdot u\, dx +
                  \int_{\Gamma} q \cdot v  \, ds\,.
\end{equation}
Further in \reff{b1} the function $g$ is given by $p$ (see \reff{a1}) through
\begin{displaymath}
g:[0,\infty)\to [0,\infty), t\mapsto g(t) = \int_0^t s\cdot p(s)\, ds ,
\end{displaymath}
where we assume that $p$ is $C^1$, $ 0\le p(t) \le p_0<\infty $,
and  $ t \mapsto t\cdot p(t)$ is 
 strictly monotonic increasing.
Then, $ 0 \le g(t)\le \frac{1}{2} p_0 \cdot t^2 $ and 
the real-valued functional 
\begin{displaymath}
G(u):=  \int_{\Omega} g (|\nabla u|)dx,  \qquad  u\in H^1(\Omega)
\end{displaymath}
is strictly convex. The Frechet derivative of $G$,
\begin{equation} \label{def-DG}
DG(u;v)=\int_{\Omega} p (|\nabla u|)  (\nabla u)^T \cdot \nabla v \, dx
\qquad  u,v\in H^1(\Omega)
\end{equation}
is Lipschitz continuous and strongly monotone in $H^1(\Omega)$ with respect to the semi-norm 
$$| v|_{H^1(\Omega)}= \|\nabla \, v \|_{L^2(\Omega)} \,, $$
that is, there exists a constant $c_G>0$ such that 
\begin{equation} \label{DG-coerc}
c_G\: \hnorm{u- v}{H^1(\Omega)}^2  \le DG(u; u-v)-DG(v; u-v)
\quad \forall  u, v \in H^1(\Omega).
\end{equation}
Analogously to \cite{CCJG-1997, MaiSte-2005}  we first define 
\begin{eqnarray*}
 \left. \begin{array}{ll}
\mathcal{L}_0:= \{ v\in H^1_{loc} (\Omega^c) : 
& \Delta v= 0 \; \mbox{in} \; H^{-1} (\Omega^c) \\[0.5ex]
& \mbox{(and for } d=2 \, \exists a \in \R
\; \mbox{such that } v \mbox{ satisfies }  (\ref{a3}))\},
 \end{array} \right. 
\end{eqnarray*}
and then the affine, hence convex set of admissible functions
\[
C := \{
(u_1,u_2)\in H^1(\Omega)\times H^1_{loc}(\Omega^c) \, : \, \;
         u_1|_{\Gamma_t} = u_2|_{\Gamma_t}+ u_0|_{\Gamma_t}  
\mbox{ and }  u_2 \in \mathcal{L}_0 \}.
\]
According to \cite[Remark 4]{CCJG-1997}, $C$ is closed in $H^1(\Omega)\times
H^1_{loc}(\Omega^c)$.
Further, we have 
\begin{eqnarray*}
 \left. \begin{array}{ll}
D\Phi ((\hat{u}_1,\hat{u}_2);(u_1,u_2))= & D G (\hat{u}_1;u_1)\, 
+ \int _{\Omega^c} \nabla \hat{u}_2 \cdot \nabla u_2 \, dx \\[0.5ex]
& - \int _{\Omega} f \cdot u_1 \, dx 
\, - \int _{\Gamma_t} q \cdot u_2|_{\Gamma_t} \, ds \,.
 \end{array} \right. 
\end{eqnarray*}

Now we can pose the HVI problem $(P_\Phi)$: Find $(\hat{u}_1, \hat{u}_2) \in C$
such that for all $(u_1, u_2) \in C$ there holds
for $\delta u_1 := u_1-\hat{u}_1 , \delta u_2 := u_2-\hat{u}_2$,  
\begin{equation} \label{HVI-1}
D \Phi ((\hat{u}_1, \hat{u}_2);( \delta u_1, \delta u_2))
 + J^0 (\gamma (\hat{u}_2- \hat{u}_1 + u_0); 
\gamma (\delta u_2 - \delta u_1)) \geq 0 \,.
\end{equation}

\begin{theorem} \label{th1}
The HVI problem $(P_\Phi)$ is equivalent to (\ref{a1})~-~(\ref{a5})
in the sense of distributions.
\end{theorem}
{\bf Proof.} First, taking into account the definition of the generalized
gradient, we note that
\begin{equation} \label{eq4}
 \int_{\Gamma_s} p(|\nabla \hat{u}_1|)\frac{\partial \hat{u}_1}{\partial n}
~ \psi \, ds \leq 
 \int _{\Gamma_s} j^0 (\cdot, (\hat{u}_2 - \hat{u}_1 + u_0)|_{\Gamma_s};
 \psi) \, ds, \,
 \forall \psi \in C^ \infty(\Gamma_s)
\end{equation}
 is the integral formulation of the nonmonotone boundary inclusion  \reff{a5}.
\\
Let $(\hat{u}_1,\hat{u}_2)\in C$ solve (\ref{HVI-1}). To show that
$(\hat{u}_1,\hat{u}_2)$ solves (\ref{a1})~-~(\ref{a5}) in the sense of
distributions, first choose $\eta \in C_0^{\infty}(\R^d)$ such that
$(u_1, u_2) := (\hat{u}_1 + \eta|_{\Omega}, \hat{u}_2+\eta|_{\Omega^c}) \in C$.
Setting $(u_1, u_2)$ in (\ref{HVI-1}) and integration by parts, implies
\begin{eqnarray*}
& 0 \,\leq & -\int_\Omega\Bigl( f +\div p(|\nabla \hat{u}_1|) \nabla
\hat{u}_1 \Bigr) \cdot \eta \, dx
  -  \int_{\Omega^c} \Delta  \hat{u}_2\cdot \eta \, dx \\[0.75ex]
& & + \ds
\int_{\Gamma}\Bigl(p(|\nabla \hat{u}_1|)\frac{\partial \hat{u}_1}{\partial n}
   -\frac{\partial \hat{u}_2}{\partial n}-q \Bigr) \cdot  \eta \, ds \,, 
  \end{eqnarray*}
since the last term in (\ref{HVI-1}) vanishes for the chosen $(u_1,u_2)$, and
moreover, $n$ pointing into $\Gamma_0 $ yields the negative sign of
$\frac{ \partial \hat{u}_2}{\partial n} $.  
Varying $\pm \eta\in C_0^{\infty} (\Omega) $ and
$\pm \eta\in C_0^{\infty} (\Omega^c) $,
shows that \reff{a1} and \reff{a2} hold in the sense of distributions. Hence,
\begin{displaymath}
0\le \int_{\Gamma}
\Bigl(  p(|\nabla \hat{u}_1|)\frac{\partial \hat{u}_1}{\partial n}
 - \frac{\partial \hat{u}_2}{\partial n}-q  \Bigr) \cdot  \eta \, ds .
\end{displaymath}
Since $\eta $ is arbitrary on $\Gamma$, we obtain 
\begin{equation} \label{eq4a}
p(|\nabla \hat{u}_1|)\frac{\partial \hat{u}_1}{\partial n} =
\frac{\partial \hat{u}_2}{\partial n}+q \qquad \mbox{a.e. on} \; \Gamma .
\end{equation}
This proves (\ref{a5_1}) and the second relation in (\ref{a4}).

Next, let $\eta_1,\eta_2\in C_0^{\infty}(\R^d) $ and consider
$(u_1, u_2):=(\hat{u}_1+ \eta_1|_\Omega, \hat{u}_2+ \eta_2|_{\Omega^c})\in C$
in an analogous way to obtain
\[
0 \leq \int_\Gamma p(|\nabla \hat{u}_1|)\frac{\partial \hat{u}_1}{\partial n}
\eta_1 -(\frac{\partial \hat{u}_2}{\partial n}+q) \eta_2 \, ds
+ \int _{\Gamma_s} j^0 (\cdot,(\hat{u}_2- \hat{u}_1 + u_0)|_{\Gamma_s}; 
(\eta_2-\eta_1)|_{\Gamma_s}) \, ds,
\]
which implies by \reff{eq4a}
\[
 \int_{\Gamma} p(|\nabla \hat{u}|)\frac{\partial \hat{u}_1}{\partial n}
 ~ (\eta_2-\eta_1) \, ds \leq 
 \int _{\Gamma_s} j^0 (\cdot,(\hat{u}_2 - \hat{u}_1 + u_0)|_{\Gamma_s}; 
(\eta_2-\eta_1)|_{\Gamma_s}) \, ds \,. 
\]
Finally, we define $\psi:=\eta_1-\eta_2$. Taking $\eta_1=\eta_2$ on $\Gamma_t$,
we have $\psi=0$ on $\Gamma_t$, but $\psi$ is arbitrary on $\Gamma_s$, what
gives (\ref{eq4}). 

Vice versa we show that (\ref{HVI-1}) follows from (\ref{a1})~-~(\ref{a5}). Let
$(u_1, u_2)$ solve (\ref{a1})~-~(\ref{a5}). Due to (\ref{a2}) and (\ref{a3}),
$(u_1, u_2) \in C$. Multiplying (\ref{a1}) and (\ref{a2}) with
differences $w_1:= v_1 - u_1$, $w_2:= v_2 - u_2$, respectively, 
where one chooses $(v_1,v_2) \in C$ arbitrarily,
and integrating by parts yields
\begin{eqnarray} 
DG(u_1;w_1) -\int_{\Gamma} p(|\nabla u_1|) 
\frac{\partial u_1}{\partial n} \, w_1 |_{\Gamma}  ds 
& = & \int_{\Omega} f \cdot w_1 \, dx, \label{eq3_1} \\[0.5ex]
\int_{\Omega^c} \nabla u_2 \cdot \nabla w_2 \, dx + \int_{\Gamma} 
\frac{\partial u_2}{\partial n} \, w_2 |_{\Gamma} ds 
& = &  0 \,. \label{eq3_2}
\end{eqnarray}
Combining (\ref{eq3_1}), (\ref{eq3_2}), 
and $p(|\nabla u_1|) \frac{\partial u_1}{\partial n} - \frac{\partial
u_2}{\partial n} = q$ on $\Gamma$,
we obtain
\[
D\Phi ((u_1,u_2); (w_1, w_2))=\int_{\Gamma} p(|\nabla u_1|) 
\frac{\partial u_1}{\partial n} \, (w_1 - w_2)|_{\Gamma} \, ds \,,
\]
where the latter integral vanishes on $\Gamma_t$ by definition of $C$.
Hence, by (\ref{a5}) (see in particular the integral  
formulation (\ref{eq4})) we conclude that for all $(v_1,v_2) \in C$,
$w_1 = v_1 - u_1$, $w_2 = v_2 - u_2$,  
\begin{eqnarray*}
& & D \Phi ((u_1, {u}_2);(w_1, w_2))
+ \int_{\Gamma_s} j^0 (\cdot,(u_2 - u_1 + u_0)|_{\Gamma_s};
(w_2 - w_1)|_{\Gamma_s})~ds \\
  &  & =  - \int_{\Gamma_s} p(|\nabla u_1|) 
\left.\frac{\partial u_1}{\partial n} (w_2 - w_1)) \right \vert_{\Gamma_s} \,
ds
\\ & & \; \; \; + \int_{\Gamma_s} j^0 (\cdot,(u_2 - u_1 + u_0)|_{\Gamma_s};
(w_2 - w_1)|_{\Gamma_s})~ds
 \\ & & \geq 0 \,,
\end{eqnarray*}
what shows that $(u_1, u_2)\in C$ solves (\ref{HVI-1}). \qed

\section{The boundary/domain HVI 
formulation of the interface problem}
\setcounter{equation}{0}

In this section we employ boundary integral operator theory
\cite{HsWe-2008,GwiSte-2018}
to rewrite  the exterior problem 
(\ref{a2})~-~(\ref{a3}) as a boundary variational inequality on $\Gamma$. 
 From now on we concentrate our analysis to the case of dimension $d =3$, since as already the distinction in the radiation condition 
(\ref{a3}) indicates, in the case $d=2$ some pecularities of boundary integral methods for exterior problems come up that need extra attention, 
 see e.g. \cite{CCJG-1997}, \cite[Sec. 12.2]{GwiSte-2018}.
As a result we arrive at an equivalent hemivariational formulation of the original interface problem \reff{a1}~-~\reff{a5} that lives on
$\Omega \times \Gamma$ and consists of a weak formulation of the nonlinear
differential operator in the bounded domain $\Omega$, the
Poincare-St\'{e}klov operator on the bounded boundary $\Gamma$, and a nonsmooth functional on the boundary part $\Gamma_s$.

To this end we need the following representation formula, see 
\cite[(1.4.5)]{HsWe-2008},\cite[(12.28)]{GwiSte-2018}.
\begin{lemma}
\label{l1}
For $u_2 \in{\cal L}_0 $ with Cauchy data $(v,\psi)$  $ $
there holds
\begin{equation}\label{c1}
u_2(x') = \frac 12 ( K v(x') - V \psi(x') ) + a, \,
x' \in\Omega^c \,, 
\end{equation}
where $V$ and $K$ denote the single layer potential and the 
double layer potential, respectively,
and $a$ is the constant appearing in \reff{a3} for $d=2$ (and
$a=0$ in (\ref{c1}) if $d > 2 $). 
\end{lemma}
Note that \reff{c1} determines $u_2$ in $\Omega^c $ as far as one knows
its Cauchy data on $\Gamma$.

Next we recall the Poincar\'{e}--Steklov operator for the exterior problem, $S:
H^{1/2}(\Gamma) \rightarrow H^{-1/2}(\Gamma) $
is a selfadjoint operator with the defining property
\begin{equation}\label{def-S}
 S(u_2|_\Gamma) = -\partial_n u_2|_\Gamma  
\end{equation}
for solutions $u_2 \in \mathcal{L}_0 $  of the Laplace equation on $\Omega^c$.
The operator $S$ enjoys the important property that it can be expressed as 
\[
S = \frac{1}{2} [ W  +
 (I-K' ) V^{-1}  (I-K) ] \,,
\]
where $I, V, K, K', W$ denote the identity, the single layer boundary integral
operator, the double layer boundary integral operator, its formal adjoint, and
the hypersingular integral operator, respectively; 
see \cite[Sec. 12.2]{GwiSte-2018} for details.

Further, $S$ gives rise to the positive definite bilinear form
$\langle S \cdot, \cdot \rangle$, 
that is, there exists a constant $c_S >0$ such that 
\begin{equation}\label{pos-def}
\langle S v, v \rangle \geq  c_S \|v\|^2_{H^{1/2}(\Gamma)} \,,
\quad  \forall v \in H^{1/2}(\Gamma) \,,
\end{equation}
where $\langle \cdot, \cdot \rangle
= \langle \cdot, \cdot \rangle_{H^{-1/2}(\Gamma) \times H^{1/2}(\Gamma) }$
extends the $L^2$ duality on $\Gamma$.

Let $E:= H^1(\Omega)\times \widetilde H^{1/2}(\Gamma_s )$
with $\widetilde H^{1/2}(\Gamma_s ) := \{ w\in H^{1/2}(\Gamma) |
\textnormal{supp}\, w\subseteq \overline\Gamma_s \} $.

We mainly follow \cite{MaiSte-2005,GMS-2021} and use the affine change of
variables
\begin{equation}\label{trafo}
(u_1,u_2) \mapsto (u,v) := (u_1 - c, u_0 + u_2|_\Gamma - u_1|_\Gamma) \in E 
\end{equation}
for a suitable $c \in \R \, (c=0,\mbox{if }d=2)$. Note that $v$ is indeed
supported in $ \overline\Gamma_s$, since the boundary condition
$(u_1 - u_2)|_{\Gamma_t} = u_0$ guarantees $v|_{\Gamma_t} = 0$.

Next, define the linear functional $\lambda \in E^*$  by 
\begin{displaymath}
\lambda(u,v) := 
 \int _{\Omega} f \cdot u \, dx
 + \langle q + Su_0, u|_{\Gamma} + v\rangle \,, \quad (u,v) \in E 
\end{displaymath}

and consider the potential energy function 
$$
{\cal E}(u,v) := G(u) + 
\frac{1}{2} \langle S(u|_\Gamma +v), u|_\Gamma  + v \rangle 
+ J(\gamma v) - \lambda(u,v) 
= \Pi(u_1, u_2 ) + C \,,
%:=  \Phi(u_1, u_2)  + J(j(u_2 - u_1 + u_0))  
$$
where $\gamma: \widetilde H^{1/2}(\Gamma_s )\rightarrow L^2(\Gamma_s )$
denotes the linear continuous embedding operator and
$ C = C(u_0,q)$ is a constant independent of $u,v$.

Then the nonsmooth, nonconvex constrained optimization problem: 
\begin{eqnarray}
\label{Energy}
\begin{array}{ll}
\mbox{minimize} & {\cal E}(u, v) \\
\mbox{subject to} & (u, v) \in E \,
\end{array}
\end{eqnarray}
leads to the following hemivariational inequality problem $(P_\mathcal{A})$:
 Find   $(\hat u,\hat v)\in E $ such that for all $(u,v) \in E$, 
\begin{equation} \label{HVIPJ}
\mathcal{A}(\hat{u}, \hat{v}; u-\hat{u}, v-\hat{v}) 
+ J^0 (\gamma \hat{v}; \gamma(v-\hat{v}))  \geq \lambda (u-\hat{u}, v-\hat{v}) \,, 
\end{equation}
where $\mathcal{A} : E \to E ^*$  is defined for all $(u,v), (u',v') \in E$ by\begin{displaymath}
\mathcal{A}(u,v)\, (u',v') = \mathcal{A}(u,v; u',v')
:= DG(u,u') + \langle S(u|_\Gamma + v),  u'|_{\Gamma} + v' \rangle \,.
\end{displaymath}

Since the nonfrictional convex smooth part of $\cal P$ coincides with the
functional $J$ of \cite{MaiSte-2005}, we immediately obtain from
\cite[Theorem 2]{MaiSte-2005} that the problems $(P_\Phi)$ and
$(P_\mathcal{A})$ are equivalent in the following way.

\begin{theorem}\label{p2}
(i) Let $(u_1,u_2)\in C $ solve $(P_\Phi)$.
Then, $(u,v)\in D$ defined by (\ref{trafo}) solves $(P_\mathcal{A})$.
\\ 
(ii) Let  $(u,v)\in D$ solve $(P_\mathcal{A})$.
%Take $a\in\R$ arbitrarily if $d=2$, whereas $a=0$ if $d = 3$. 
Define $u_1:=u$ and $u_2$ by the representation formula (\ref{c1})
with $(u|_{\Gamma}+v-u_0,-S(u|_{\Gamma}+v-u_0))$ replacing $(v,\psi)$, i.e. 
$u_2:=\frac{1}{2} \left( K(u|_{\Gamma} + v-u_0)
+ V(S(u|_{\Gamma}-v+u_0 ))\right)$.
Then, $(u_1,u_2)\in C$ solves $(P_\Phi)$.
\end{theorem}

Thanks to the strong monotonicity of the nonlinear operator $DG$ in
$H^1(\Omega)$ with respect to the semi-norm $| \cdot|_{H^1(\Omega)}= \|\nabla
\, \cdot \|_{L^2(\Omega)}$ and the positive definiteness of the
Poincar\'{e}--Steklov operator $S$ the following strong monotonicity 
property can be derived, see \cite[Lemma 4.1]{CCJG-1997}. 

\begin{lemma} \label{l-coerc}
There exists a constant $c_{\cal A} >0$
such that for all
$v,v' \in \tilde H^{1/2}(\Gamma_s)$ and all $u,u' \in H^1(\Omega)$
there holds
\begin{eqnarray*}
\lefteqn{ c_{\cal A} \cdot \Vert (u-u',v-v')
\Vert^2_{ H^1(\Omega)\times \tilde H^{1/2}(\Gamma_s) }}\\
&\le &  DG(u; u-u')- DG(u';u-u') \\ &+&
 \langle S(u|_\Gamma + v - u'|_\Gamma - v')}{u|,
u_\Gamma + v - u'|_\Gamma - v' \rangle.
\end{eqnarray*}
\end{lemma}

Moreover, following the arguments in \cite{OvGw-Rassias-14} and
using the compact embedding 
$\widetilde H^{1/2}(\Gamma_s ) \subset L^2 (\Gamma_s)$
 it can be easily seen that the  functional 
$\varphi : \widetilde H^{1/2}(\Gamma_s ) \times \widetilde H^{1/2}(\Gamma_s )
 \to \R$ defined by
\[
 \varphi (v,\tilde{v}) 
:= \int _{\Gamma_s} j^0(\cdot, v;  \tilde{v} - v)\, ds
\]
is pseudomonotone and upper semicontinuous.   
Thus the concrete hemivariational inequality problem $(P_\mathcal{A})$
is covered by the general theory exhibited in Section \ref{inter};
existence and under the smallness condition uniqueness hold for
$(P_\mathcal{A})$.

\section{Extended real-valued HVIs - existence, uniqueness, and stability} 
\label{extend}
\setcounter{equation}{0}

In view of the subsequent study of optimal control problems  in Section \ref{optconinter} governed by  the interface problem which we have investigated in the previous sections, we broaden the scope of analysis and consider extended real-valued HVIs: 
Find $\hat v\in \textnormal{dom } F $ such that
\begin{equation}\label{a66}
{\cal A}(\hat v) (v- \hat v) 
+ J^0(\gamma \hat v; \gamma v- \gamma \hat v)  %\varphi(\hat{v},v)
%\ge \lambda(v-\hat v)
+ F(v) -F(\hat v)  \ge 0
\qquad \forall v \in V\,.
\end{equation}
Here $V$ is a real reflexive Banach space,
the nonlinear operator ${\cal A} : V\to V^* $ is a monotone operator,
  $\gamma := \gamma_{V \to X}$ 
with  $X$ a real Hilbert space (in the interface problem 
we have $X= L^2 (\Gamma_s)$) denotes a linear continuous operator,
$J^0$ stands for the generalized directional derivative of
a real-valued locally Lipschitz functional $J$,
and now in addition, $F: V \rightarrow \R \cup \{+\infty \}$ is a convex lower semicontinuous function that is supposed to be
proper (i.e. $F \not\equiv \infty$ on ${\cal C}$).
This means that the effective domain of $F$
in the sense of convex analysis (\cite{RTR}),
$$
\textnormal{dom } F := \{v \in V: F(v) < + \infty  \}
$$
 is nonempty, closed and convex.
To resume the HVI (\ref{a6}) of  Sect. \ref{inter}, 
let $F(v) := \lambda(v) + \chi_{\cal C} (v)$,
where $\lambda \in V^*$ and
\[
\chi_{\cal C} (v) :=
~~\left\{ \begin{array}{ll}
  0  & \hbox{ if } v \in {\cal C} \\
  + \infty & \hbox{ elsewhere } 
  \end{array} \right.
\] 
 is the indicator function on ${\cal C}$
in the sense of convex analysis (\cite{RTR}).

Next similarly to (\ref{bif1}) in Sect. \ref{inter} define
\begin{equation}\label{bif2} 
\varphi(v,w) := 
{\cal A}(v) (w- v) + J^0(\gamma v; \gamma w - \gamma v)  
\end{equation}
and apply Proposition \ref{prop1}.
Thus under the assumptions \reff{cc1} and \reff{cc2}, the above HVI 
(\ref{a66}) falls into the framework of an 
{\it extended real-valued  equilibrium problem of monotone type} in the sense of \cite{Gwi-22}. Clearly  strong monotonicity implies uniqueness. 
Note by the separation theorem it can be shown that any 
convex proper lower semicontinuous function
$\phi: V \rightarrow \R \cup \{+\infty \}$ 
is conically minorized, that is, it enjoys the estimate
$$
\phi(v) \ge - c_\phi(1+ \|v\|)  
$$
with some $c_\phi > 0$.  
Hence strong monotonicity implies the asymptotic coercivity 
condition in  \cite{Gwi-22}, too. 
Thus the existence result \cite[Theorem 5.9]{Gwi-22} applies to the HVI (\ref{a66}) to conclude the following

\begin{theorem} \label{theo-11}
 Suppose \reff{cc1} and \reff{cc2}.
Then the  HVI (\ref{a66}) is uniquely solvable.
\end{theorem}

By this solvability result, we can introduce the solution map
${\cal S}$ by ${\cal S}(F) := \hat v$,
the solution of (\ref{a66}). Next we investigate the stability 
of the solution map ${\cal S}$ with respect to the extendend
real-valued function $F$. Here we follow the concept of 
epi-convergence in the sense of Mosco  \cite{Mos1969,AttBut-2014}
  ("Mosco convergence"). 
Let  $F_n ~ (n \in \N), F: V \rightarrow \R \cup \{+\infty \}$
be convex lower semicontinuous proper functions. Then
$F_n$ are called to converge to $F$ in the Mosco sense,
written $F_n \stackrel{\rm M}{\longrightarrow} F$, if and only if
the subsequent two hypotheses hold:
\begin{enumerate}
\item[{\bf (M1)}] If $ v_n \in V ~(n \in \N)$
weakly converge to $v$ for $n \to \infty$, then
$$
F(v) \le \liminf_{n \to \infty} F_n (v_n) \,.
$$
\item[{\bf (M2)}] For any $v \in  V$ there exist  $ v_n \in  V~(n \in \N)$ 
 strongly converging to $v$ for $n \to \infty$ such that
$$
F(v) = \lim_{n \to \infty} F_n (v_n) \,.
$$
\end{enumerate}
In view of our later applications it is not hard to require that the functions
$F_n$ are uniformly conically minorized, that is, there holds the  estimate
\begin{equation}\label{minor}
F_n(v) \ge - d_0 (1+ \|v\|), ~\forall n \in \N, v \in V
\end{equation}  
with some $d_0 \ge 0$. Moreover similar to \cite{Sof2018},
in addition to the one-sided Lipschitz continuity (\ref{cc1}),
we assume that the  locally Lipschitz function $J$ satisfies the following
growth condition
\begin{equation}\label{cc3}
\| \zeta\|_{X^*} \le d_J (1 +  \|z\|_X) 
\qquad \forall z \in X, \, \zeta \in \partial J(z)
\end{equation}
for some $d_J > 0$, what is immediate from the growth condition
(\ref{as-j}) for the integrand $j$.

Now we are in the position to state the main result of this section which extends the stability result of \cite{Gwi-95} for monotone variational inequalities to  extended real-valued HVIs 
with an unperturbed bifunction $\varphi$ in the coercive situation.

\begin{theorem} \label{STABLE}
Suppose that the operator
${\cal A}$ is continuous and strongly monotone with monotonicity constant $c_{\cal A} > 0$, the linear operator  $\gamma$ is compact, the
generalized directional derivative $J^0$
satisfies the  one-sided Lipschitz condition (\ref{cc1})
and the growth condition (\ref{cc3}).
Moreover, suppose the smallness condition (\ref{cc2}).
Let  $F, F_n: V \rightarrow \R \cup \{+\infty \}~ (n \in \N)$ be convex lower semicontinuous proper functions 
that satisfy the lower estimate(\ref{minor}); 
let $F_n \stackrel{\rm M}{\longrightarrow} F$.
Then strong convergence ${\cal S}(F_n) \rightarrow {\cal S}(F)$ holds. 
\end{theorem}
{\bf Proof.\/} We divide the proof in three parts. We first show 
that the $\hat u_n = {\cal S}(F_n)$ are bounded, before we can establish the convergence result. 
In the following $c_0,c_1,\ldots \,$ are generic positive constants.

{\bf (1)} {\em The sequence $\{ \hat u_n\}\subset V$ is bounded.\/}.

By definition, $\hat u_n$ satisfies for all $v \in V$,
\begin{equation}\label{a66n}
{\cal A}(\hat u_n) (v- \hat u_n) 
+ J^0(\gamma \hat u_n; \gamma v- \gamma \hat u_n)  
+ F_n(v) -F_n(\hat u_n)  \ge 0 \,.
\end{equation}
Now let $v_0$ be an arbitrary element of $\textnormal{dom } F$.
Then by Mosco convergence, (M2),
 there exist  $ v_n \in \textnormal{dom } F_n~(n \in \N)$ 
such that for $n \to \infty$  the strong convergences hold
 \begin{equation} \label{ass-M1} 
v_n \rightarrow v_0; \, F_n (v_n) \rightarrow F(v_0)  \,.
\end{equation} 
Let $n \in \N$. Then insert $v = v_n$ in (\ref{a66n}) and obtain
$$ 
 {\cal A}(\hat u_n) (\hat u_n - v_n) 
 \le J^0(\gamma \hat u_n; \gamma v_n - \gamma \hat u_n)   
+ F_n(v_n) - F_n(\hat u_n) \,.
$$
Write ${\cal A}(\hat u_n) = 
{\cal A}(\hat u_n) - {\cal A}(v_n) + {\cal A}(v_n) $
and use the strong monotonicity of the operator ${\cal A}$ 
and the estimate (\ref{minor}) to get
\begin{eqnarray}
\label{est-1} 
&& c_{\cal A} \|\hat u_n - v_n \|^2_V \\[0.5ex] \nonumber 
&& \le \| {\cal A}(v_n) \|_{V^*} \|\hat u_n - v_n \|_V
+ F_n(v_n) + d_0 (1+ \| \hat u_n \|) \\[0.5ex] \nonumber 
&& \quad  + ~ J^0(\gamma \hat u_n; \gamma v_n - \gamma \hat u_n) \,. 
  \end{eqnarray}
On the other hand write 
\begin{eqnarray*}
&& J^0(\gamma \hat u_n; \gamma v_n - \gamma \hat u_n) \\[0.5ex]
&& = J^0(\gamma \hat u_n; \gamma v_n - \gamma \hat u_n) +
J^0(\gamma v_n; \gamma \hat u_n - \gamma v_n)%  \\[0.5ex] && \quad
 - J^0(\gamma v_n; \gamma \hat u_n - \gamma v_n) \,. 
 \end{eqnarray*}
Hence by the  one-sided Lipschitz condition (\ref{cc1}),
\begin{equation}\label{est-2}
J^0(\gamma \hat u_n; \gamma v_n - \gamma \hat u_n)
 \le c_J \| \gamma \|^2 \| v_n - \hat u_n \|_V^2 
- J^0(\gamma v_n; \gamma \hat u_n - \gamma v_n) \,.
\end{equation}
Further by (\ref{cc3}),
\begin{eqnarray} \nonumber
&& - J^0(\gamma v_n; \gamma \hat u_n - \gamma v_n)
\le \max_{\zeta \in \partial J(\gamma v_n)} 
\| \zeta\|_{X^*} \|\gamma \hat u_n - \gamma v_n\|_X \\[0.5ex] \label{est-3}
&& \le d_J \| \gamma \| (1+ \| \gamma \|  \|v_n\|_V) \|\hat u_n - v_n\|_V \,.  
\end{eqnarray}
By the convergences (\ref{ass-M1}), $|F_n(u_n)| \le c_0 ,  \| {\cal A}(v_n) \|_{V^*} \le c_1, 
 \|v_n\|_V \le c_2$. Thus (\ref{est-1}), (\ref{est-2}), (\ref{est-3})
result in
$$
( c_{\cal A} -  c_J \| \gamma \|^2) \|\hat u_n - v_n \|^2_V \le
c_0 + [c_1 +  d_J \| \gamma \| (1+ c_2 \| \gamma \|)] \|\hat u_n - v_n \|_V
+ d_0 (1+ \| \hat u_n \|) \,.
$$
Hence by the smallness condition (\ref{cc2}), a contradiction argument
proves the claimed boundedness of $\{ \hat u_n\}$.

{\bf (2)} {\em  $ \hat u_n ={\cal S} (F_n)$ converges weakly to
$ \hat u ={\cal S} (F)$ for $n \to \infty$ .\/}.

To prove this claim we employ a "Minty trick" similar to the proof of 
\cite[Prop.3.2]{Gwi-22} using the monotonicity of the operator
$\cal A$.

Take $v \in V$ arbitrarily. By $(M2)$ there exist  $ v_n \in  V~(n \in \N)$
such that  
\begin{equation}\label{est-4}
 \lim_{n \to \infty} v_n = v; \, \lim_{n \to \infty} F_n (v_n) = F(v)  
\end{equation}
We test the inequality (\ref{a66n}) with $v_n$, use the monotonicity of the operator $\cal A$, and obtain 
\begin{equation}\label{est-5}
{\cal A}(v_n) (v_n - \hat u_n) 
+ J^0(\gamma \hat u_n; \gamma v_n - \gamma \hat u_n) \ge 
F_n(\hat u_n) -  F_n(v_n)  \,.
\end{equation}
On the other hand, by the previous step, there exists a subsequence 
$\{ \hat u_{n_k} \}_{k \in \N}$ that converges weakly to some 
$\tilde u  \in \textnormal{dom } F\subset V$.
Further since $\gamma$ is completely continuous, 
$\gamma \hat u_{n_k} \rightarrow \gamma \tilde u$.
Thus the continuity of $\cal A$, the upper semicontinuity
of $(y,z) \in  X \times X \mapsto J^0(y;z)$, $(M1)$, and (\ref{est-4}) 
entail together with (\ref{est-5})
\begin{eqnarray*}
&& {\cal A}(v) (v - \tilde u) 
+ J^0(\gamma \tilde u; \gamma v - \gamma \tilde u)\\[0.5ex]
&& \ge \lim_{k \to \infty} {\cal A}(v_{n_k}) (v_{n_k} - \hat u_{n_k}) 
+  \limsup_{k \to \infty} 
J^0(\gamma \hat u_{n_k};\gamma v_{n_k} - \gamma \hat u_{n_k}) \\[0.5ex]
&& \ge \liminf_{k \to \infty} F_n(\hat u_n) - \lim_{k \to \infty} F_n(v_n) \\[0.5ex]
&& \ge F(\tilde u) - F(v) \,. 
 \end{eqnarray*}
Hence for $v \in  \textnormal{dom } F$ fixed, for arbitrary 
$s \in [0,1)$ and $w_s := v + s(\tilde u -v) \in \textnormal{dom } F$
inserted above, 
the positive homogeneity of $J^0(\gamma \tilde u; \cdot)$ and the
convexity of $F$ imply after division by the factor $(1-s) > 0$ 
$$
{\cal A}(w_s) (v - \tilde u) 
+ J^0(\gamma \tilde u; \gamma v - \gamma \tilde u) + F(v) \ge  F(\tilde u) \,.
$$
Letting $s \rightarrow 1$, hence $w_s \rightarrow \tilde u$,
${\cal A}(w_s) \rightarrow {\cal A}(\tilde u)$ results in 
$$
{\cal A}(\tilde u) (v - \tilde u) 
+ J^0(\gamma \tilde u; \gamma v - \gamma \tilde u) + F(v) \ge  F(\tilde u) 
\quad \forall v \in \textnormal{dom } F \,.
$$
This shows by uniqueness that $\tilde u = {\cal S} (F)$ and the entire
sequence $\{ \hat u_n\}$ converges weakly to $\hat u = {\cal S} (F)$. 

{\bf (3)} {\em  $ \hat u_n ={\cal S} (F_n)$ converges strongly to
$ \hat u ={\cal S} (F)$ for $n \to \infty$.\/}

By $(M2)$ there exist  $ u_n \in  V~(n \in \N)$
such that  
\begin{equation}\label{est-6}
(i) \, \lim_{n \to \infty} u_n = \hat u; 
(ii)\, \lim_{n \to \infty} F_n (u_n) = F(\hat u) \,.  
\end{equation}
Test the inequality (\ref{a66n}) with $u_n$, use the strong monotonicity of the operator $\cal A$, and obtain 
\begin{eqnarray} \nonumber
&& {\cal A}(u_n) (u_n - \hat u_n) 
+ J^0(\gamma \hat u_n; \gamma u_n - \gamma \hat u_n) \\[0.5ex] \label{est-7}
&& + F_n(u_n) - F_n(\hat u_n) \ge  c_{\cal A} \, \| u_n - \hat u_n \|^2 \,.  
\end{eqnarray}
Analyze the summands in  (\ref{est-7}) separately:
By (\ref{est-6}) (i), ${\cal A}(u_n) \rightarrow {\cal A}(\hat u)$, hence
$$
 \lim_{n \to \infty} {\cal A}(u_n) (u_n - \hat u_n) = 0 \,.
$$
By the upper semicontinuity
of $(y,z) \in  X \times X \mapsto J^0(y;z)$, $(M1)$
and by the complete continuity of $\gamma$,
$$
\limsup_{n \to \infty} 
 J^0(\gamma \hat u_n; \gamma u_n - \gamma \hat u_n) \le 0 \,.
$$
By (\ref{est-6}) (ii) and by $M1$,
$$
\limsup_{n \to \infty} 
[ F_n(u_n) - F_n(\hat u_n) ] \le 0 \,.
$$
Thus  from (\ref{est-7}) finally by the triangle inequality,
$$
0 \le \| \hat u_n - \hat u \| \le
  \| \hat u_n - u_n \|  +  \| u_n - \hat u \| \rightarrow 0
$$
and the theorem is proved. 
\qed

%\begin{remark} \label{rem-MS}
To conclude this section let us compare the above Theorem \ref{STABLE} with a similar stability result  of \cite[Theorem 6]{Sof2018}.   
There one has the special setting of
$F(v) :=  \chi_{\cal C} (v) + \lambda(v)  $,
where ${\cal C} \subset V$  is closed convex and $\lambda \in V^*$
is given by $\lambda := \kappa^*f$ with $f \in X^*$, $\kappa^*$ the adjoint  
to the linear operator $\kappa: V \rightarrow X$, which is assumed to be 
completely continuous or equivalently compact, 
see \cite[(4.3)]{Sof2018}.

Likewise for $ n \in \N$, one has convex closed sets
${\cal C}_n \subset V$ and linear forms 
$\lambda_n := \kappa^*f_n$ with $f_n \in X^*$ giving rise to 
$F_n(v) :=  \chi_{{\cal C}_n} (v) + \lambda_n(v)  $. 

Note by the Schauder theorem, see e.g. \cite[3.7.17]{PaWi2018},
the adjoint $\kappa^*: X^* \rightarrow V^* $ is compact.
Hence the assumed weak convergence $f_n \rha f$ in $X^*$
entails the strong convergence $ \lambda_n  \rightarrow \lambda $ in $V^*$,
thus further the lower estimate (\ref{minor}),
in view of $\chi_{{\cal C}_n} \ge 0$.
Moreover, the Mosco convergence
$F_n \stackrel{\rm M}{\longrightarrow} F$ follows at once from the
assumed   Mosco convergence 
${\cal C}_n \stackrel{\rm M}{\longrightarrow} {\cal C}$,
namely from the hypotheses
\begin{enumerate}
\item[{\bf (m1)}] If $ v_n \in {\cal C}_n ~(n \in \N)$
weakly converge to $v$ for $n \to \infty$, then 
$ v \in {\cal C}$.
\item[{\bf (m2)}] For any $v \in {\cal C}$ there exist 
$ v_n \in {\cal C}_n ~(n \in \N)$
strongly converging to $v$ for $n \to \infty$. 
\end{enumerate}
% \\[1.5ex]
On the other hand an inspection of the  above proof of Theorem \ref{STABLE} shows that is enough to demand
for $ {\tilde J}^0(v;w) = J^0(\gamma v; \gamma w)$  
 the generalized directional derivative of
the  real-valued locally Lipschitz functional $\tilde J(v):= J(\gamma v)$
that
\begin{eqnarray*}
&& u_n \rha u \textnormal{ in } V \textnormal{ and }
  v_n \rightarrow v \textnormal{ in } V \\[0.5ex]
&& \Rightarrow \limsup  {\tilde J}^0(u_n;v_n - u_n) \le {\tilde J}^0(u;v-u) \,. 
 \end{eqnarray*}
This abstract condition \cite[(4.2)]{Sof2018} is derived in the above proof of Theorem \ref{STABLE} from the compactness of $\gamma$ and
the upper semicontinuity of $(y,z) \in  X \times X \mapsto J^0(y;z)$.\\[1.5ex] 
Concerning the monotone operator ${\cal A} : V\to V^* $,
we only require its norm continuity, not needing Lipschitz continuity.
More importantly, we can also dispense with the  condition \cite[(4.1)]{Sof2018}:
 \begin{eqnarray*}
&& u_n \rha u \textnormal{ in } V \textnormal{ and }
  v_n \rightarrow v \textnormal{ in } V \\[0.5ex]
&& \Rightarrow \limsup \lan  {\cal A} u_n, u_n - v_n \ran
 \ge \lan {\cal A} u, u-v \ran \,. 
 \end{eqnarray*}
It seems that this condition forces an elliptic
operator, 
which stems from an elliptic pde on the domain, 
to be linear.

%\end{remark}

\section{Some optimal control problems governed by the interface problem} 
\label{optconinter}
\setcounter{equation}{0}
In this section we  rely heavily on the stability result
Theorem \ref{STABLE} 
and present an unified approach to existence results 
 for various optimal control problems governed 
by the interface problem which was described
and studied in sections 2 -- 4.
For convenience, let us recall the boundary/domain HVI 
formulation $(P_\mathcal{A})$ of the interface problem:
 Find   $(\hat u,\hat v)\in E $ such that for all $(u,v) \in E$, 
\begin{equation} \label{HVIPJ-2}
\mathcal{A}(\hat{u}, \hat{v}; u-\hat{u}, v-\hat{v}) 
+ J^0 (\gamma \hat{v}; \gamma(v-\hat{v}))  \geq \lambda (u-\hat{u}, v-\hat{v}) \,. 
\end{equation}
Here $E= H^1(\Omega)\times \widetilde H^{1/2}(\Gamma_s )$
with $\widetilde H^{1/2}(\Gamma_s ) = \{ w\in H^{1/2}(\Gamma) |
\textnormal{supp}\, w\subseteq \overline\Gamma_s \} $
on the bounded domain $\Omega$ and the boundary part $\Gamma_s$.
The operator $\cal A$ is given for all $(u,v), (u',v') \in E$ by\begin{displaymath}
\mathcal{A}(u,v)\, (u',v') = \mathcal{A}(u,v; u',v')
= DG(u,u') + \langle S(u|_\Gamma + v),  u'|_{\Gamma} + v' \rangle \,,
\end{displaymath}
see (\ref{def-DG}),  (\ref{def-S}).
$ J^0 $ denotes the generalized directional derivative of
the Lipschitz integral function $J$, 
see (\ref{def-J0}), (\ref{def-J}) stemming from
 the generally nonmonotone, set-valued transmission condition (\ref{a5}).
$\gamma: \widetilde H^{1/2}(\Gamma_s )\rightarrow L^2(\Gamma_s )$
denotes the linear continuous embedding operator which is compact.
The linear functional $\lambda \in E^*$  is defined 
for  $(u,v) \in E$ by 
\begin{displaymath}
\lambda(u,v) = 
 \int _{\Omega} f \cdot u \, dx
 + \langle q + Su_0, u|_{\Gamma} + v\rangle \,, 
\end{displaymath}
where $\langle \cdot, \cdot \rangle
= \langle \cdot, \cdot \rangle_{H^{-1/2}(\Gamma) \times H^{1/2}(\Gamma) }$
extends the $L^2$ duality on $\Gamma$.

Now for simplicity we set $u_0 := 0$ and impose for the data 
$f, q$ that $f \in L^2(\Omega)$ and $q \in L^2(\Gamma)$. Thus
we can write 
\begin{equation} \label{linform}
\lambda(u,v) = 
 \langle f, \kappa u \rangle_{L^2(\Omega) \times L^2(\Omega)}
 + \langle q, \tau u + \iota v \rangle_{L^2(\Gamma) \times L^2(\Gamma)} \,, 
\end{equation}
where $\kappa: H^1(\Omega) \to L^2(\Omega) $,
$\iota: H^{1/2}(\Gamma) \to L^2(\Gamma) $ are
 linear compact embedding operators and 
$\tau: H^1(\Omega) \to L^2(\Gamma)$ is a
linear compact trace operator.

In the following three optimal control problems we stick to a misfit
 functional with a given target $(u_d,v_d)  \in E$ 
as the simplest case of a cost functional and regularize this functional
by the norm of the control with a given regularization parameter $\rho > 0$. The subsequent analysis can be extended
to cover more general cost functionals under appropriate
lower semicontinuity and coerciveness assumptions, 
see e.g. \cite[sec. 5]{Sof2018} , see also 
\cite[sec. 4.4]{Tr2010}, 
and also to a more general setting of regularization, 
see e.g. \cite[II Sect. 7.5 (7.51)]{Li1971}, \cite[(26)]{GJKS18}.

\subsection{A distributed optimal control problem 
governed by the interface problem}

Here we control by $f \in L^2(\Omega)$ distributed on the domain $\Omega$.
Thus in the abstract setting of section 5,  we choose the 
convex functional $F$ as the linear functional
$$
F(u,v) :=  \langle f, \kappa u \rangle_{L^2(\Omega) \times L^2(\Omega)}
=  (\kappa^* f) (u), \, F = (\kappa^* f,0) \in E^* \,.
$$
By the abstract existence and uniqueness result of Theorem \ref{theo-11}
 we have the control-to-state map $f \in L^2(\Omega) \mapsto
{\cal S}(f) := (\hat u,\hat v)$,
the solution of (\ref{HVIPJ-2}). Thus we can pose the optimal 
control problem  $(\mathcal{OCP})_1$:
\begin{eqnarray*}
%\label{ocp-1}
\begin{array}{ll}
\mbox{minimize } & I_1(f) := 
\frac{1}{2} \, \| {\cal S}(f) - (u_d,v_d) \|_E^2  
+ \frac{\rho}{2} \, \|f\|^2_{L^2(\Omega)} \\[0.75ex]
\mbox{subject to } &  f \in L^2(\Omega) \,,
\end{array}
\end{eqnarray*}
for which we can prove the following existence result.

\begin{theorem} \label{OCP-1}
Suppose that the generalized directional derivative $J^0$
satisfies the  one-sided Lipschitz condition (\ref{cc1})
and the growth condition (\ref{cc3}).
Moreover, suppose the smallness condition (\ref{cc2})
with the monotonicity constant $c_{\cal A}$ from Lemma \ref{l-coerc}.
Then there exists an optimal control 
to the optimal control problem  $(\mathcal{OCP})_1$.
\end{theorem}
{\bf Proof.\/} The proof follows the standard pattern of existence
proofs in optimal control.

Since the cost function  $I_1$ is  bounded below, 
$$
I_1^* := \inf (\mathcal{OCP})_1 \in \R \,.
$$ 
Let $\{f_n\}_{n\in\N}$ be a minimizing sequence of  $(\mathcal{OCP})_1$, e.g.
construct $f_n$ by $I_1(f_n) < I_1^* + \frac{1}{n}$.
Since the misfit term in $I_1$ is nonnegative, 
$\|f_n\|^2_{L^2(\Omega)}$ is bounded. By reflexivity of $L^2(\Omega)$, 
we can extract a subsequence of $\{f_n\}_{n\in\N}$
 also denoted $\{f_n\}_{n\in\N}$ that weakly converges to some 
$\hat f \in L^2(\Omega)$. 
Since $\kappa^*$ is completely continuous, we have strong convergence
$F_n := (\kappa^* f_n,0) \rightarrow \hat F := (\kappa^* \hat f,0)$  in $E^*$.
Thus the linear continuous functionals $F_n$ satisfy 
the lower estimate (\ref{minor}).
To show {\bf (M1)}, let $(u_n,v_n) \rha (u,v)$ in $E$ for $n \to \infty$.
Then clearly, $\hat F(u,v) =  \lim_{n \to \infty} F_n (u_n,v_n)$.
To show {\bf (M2)}, choose for any $(u,v) \in E$, simply
$(u_n,v_n) := (u,v) \in E$. Then clearly $(u_n,v_n) \rightarrow (u,v)$ 
and $ \hat F(u,v) = \lim_{n \to \infty} F_n (u_n,v_n)$. Hence
$F_n \stackrel{\rm M}{\longrightarrow} \hat F$ and 
Theorem \ref{STABLE} applies; it yields
${\cal S}(f_n) \rightarrow {\cal S}(\hat f)$. Thus in view of the weak 
lower semicontinuity of the norm, 
$$I_1^* \le I_1(\hat f) \le \liminf_{n \to \infty} I_1(f_n) \le I_1^* 
$$
and the theorem is proved.  
\qed

\subsection{A boundary optimal control problem and a
simultaneous distributed-boundary optimal control problem
governed by the interface problem}

Now we control by $q \in L^2(\Gamma)$ on the boundary $\Gamma$.
Thus in the abstract setting of section 5,  we now choose the 
convex functional $F$ as the linear functional
$$
F(u,v) :=  \langle q, \tau u  + \iota v \rangle_{L^2(\Gamma) \times L^2(\Gamma)}
=  (\tau^* q, \iota^* q) (u,v), \, F = (\tau^* q, \iota^* q) \in E^* \,.
$$
By the abstract existence and uniqueness result of Theorem \ref{theo-11}
 we have the control-to-state map $q \in L^2(\Gamma) \mapsto
{\cal S}(q) := (\hat u,\hat v)$,
the solution of (\ref{HVIPJ-2}). Thus we can pose the optimal 
control problem  $(\mathcal{OCP})_2$:
\begin{eqnarray*}
%\label{ocp-2}
\begin{array}{ll}
\mbox{minimize } & I_2(q) := 
\frac{1}{2} \, \| {\cal S}(q) - (u_d,v_d) \|_E^2  
+ \frac{\rho}{2} \,  \|q\|^2_{L^2(\Gamma)} \\[0.75ex]
\mbox{subject to } &  q \in L^2(\Gamma) \,,
\end{array}
\end{eqnarray*}
for which we can prove the following existence result.

\begin{theorem} \label{OCP-2}
Suppose that the generalized directional derivative $J^0$
satisfies the  one-sided Lipschitz condition (\ref{cc1})
and the growth condition (\ref{cc3}).
Moreover, suppose the smallness condition (\ref{cc2})
with the monotonicity constant $c_{\cal A}$ from Lemma \ref{l-coerc}.
Then there exists an optimal control 
to the optimal control problem  $(\mathcal{OCP})_2$.
\end{theorem}
{\bf Proof.\/} The proof follows from arguments similar to those that were given in the proof of Theorem \ref{OCP-1}. So the details are omitted.
\qed

Let us remark that we can also treat the 
simultaneous distributed-boundary optimal control problem
$(\mathcal{OCP})_3$ as in \cite{CHNZ2022}, 
here driven by the interface problem:

\begin{eqnarray*}
%\label{ocp-3}
\begin{array}{ll}
\mbox{minimize } & I_3(f,q) := 
\frac{1}{2} \,  \| {\cal S}(f,q) - (u_d,v_d) \|_E^2  
+ \frac{\rho}{2} \, [ \|f\|^2_{L^2(\Omega)} + \|q\|^2_{L^2(\Gamma)} ] \\[0.75ex]
\mbox{subject to } & f \in L^2(\Omega), \, q \in L^2(\Gamma) \,,
\end{array}
\end{eqnarray*}
where now we have the control-to-state map 
$(f,q) \in L^2(\Omega) \times L^2(\Gamma) \mapsto
{\cal S}(f,q) := (\hat u,\hat v)$,
the solution of (\ref{HVIPJ-2}). 
By analogous reasoning we obtain an optimal solution to
$(\mathcal{OCP})_3$. The details are omitted.
 
\subsection{An optimal control problem driven by a
bilateral obstacle interface problem} 

To conclude this section we investigate a related bilateral obstacle interface problem and the associated optimal control problem similar to \cite{LYZe2021}. 
First for the strong formulation, we modify in the interior part 
$\Omega \subset \R^3$ the nonlinear partial differential equation (\ref{a1})
to the obstacle problem: Find $u=u(x) \in [\ul u(x), \ol u(x)]$ such that
\begin{eqnarray}
\label{a-4-1} \left.
\begin{array}{lll}
- \div \Bigl( p(| \nabla u |)\cdot \nabla u \Bigr) \ge f &
\mbox{ if }  u = \ul u & \mbox{ a.e. in } \Omega, \\[0.75ex]
- \div \Bigl( p(| \nabla u |)\cdot \nabla u \Bigr) = f &
\mbox{ if }  \ul u < u < \ol u & \mbox{ a.e. in } \Omega, \\[0.75ex]
- \div \Bigl( p(| \nabla u |)\cdot \nabla u \Bigr) \le f &
\mbox{ if }  u = \ol u & \mbox{ a.e. in } \Omega, 
\end{array} \right\}
\end{eqnarray}
where the obstacle functions $\ul u, \ol u\in H^1(\Omega)$ with  
$\ul u \le \ol u\mbox{ a.e. in } \Omega$ are given. In the exterior part  $\Omega^c$, we consider still the Laplace equation (\ref{a2})
with the radiation  condition (\ref{a3}).
The transmission conditions (\ref{a4}), (\ref{a5_1}), and (\ref{a5})
remain in force. By the variational analysis in sections (3) and (4)
using boundary integral methods we arrive at
the following hemivariational inequality problem 
$(P_\mathcal{A,C})$:
 Find   $(\hat u,\hat v)\in {\cal C} $ 
such that for all $(u,v) \in {\cal C} $,  
\begin{equation} \label{HVI-obs}
\mathcal{A}(\hat{u}, \hat{v}; u-\hat{u}, v-\hat{v}) 
+ J^0 (\gamma \hat{v}; \gamma(v-\hat{v}))  \geq \lambda (u-\hat{u}, v-\hat{v}) \,, 
\end{equation}
where the operator $\mathcal{A}$, 
the generalized directional derivative $J^0$,
the linear continuous embedding operator $\gamma$,
the  linear functional $\lambda$  are defined as before and above 
in this section, and where now the constraint set
\begin{equation}
{\cal C} := {\cal C}_{\ul u,\ol u} :=
\{(u,v) \in E | \, \ul u \le u \le \ol u  \mbox{ a.e. in } \Omega  \}    
\label{a-4-2} 
\end{equation}
is closed and convex. This gives rise to the closed convex functional
%\begin{equation}
$$
F := F_{\ul u,\ol u} := \chi_{\cal C} = \chi_{{\cal C}_{\ul u,\ol u}} \,. 
$$ %\label{a-4-3}  \end{equation}

Here we control by the obstacles $\ul u, \ol u$  distributed on the domain $\Omega$ where as in \cite{LYZe2021} we impose the regularity 
$\ul u, \ol u \in H^2(\Omega)$ 
 and introduce the admissible set
$$
U_{\mbox \small {ad}} := \{(\ul u, \ol u)|\,\ul u \le \ol u \mbox{ a.e. in } \Omega\} \,.
$$
By the abstract existence and uniqueness result of Theorem \ref{theo-11}
 we have the control-to-state map $(\ul u, \ol u) \in  U_{\mbox \small {ad}} \mapsto {\cal S}(\ul u, \ol u) := (\hat u,\hat v)$,
the solution of (\ref{HVI-obs}). Thus we can pose the optimal 
control problem  $(\mathcal{OCP})_4$:
\begin{eqnarray*}
\begin{array}{ll}
\mbox{minimize } & I_4(\ul u, \ol u) := 
\frac{1}{2} \, \| {\cal S}(\ul u, \ol u) - (u_d,v_d) \|_E^2  
+ \frac{\rho}{2} \, [\|\ul u \|^2_{H^2(\Omega)} + \|\ol u\|^2_{H^2(\Omega)}]
 \\[0.75ex] \mbox{subject to } &  (\ul u, \ol u) \in U_{\mbox \small {ad}} \,.
\end{array}
\end{eqnarray*}

An essential ingredient in the subsequent proof of the existence of an optimal control to $(\mathcal{OCP})_4$ is the Mosco convergence of constraint sets.
For that latter result we exploit  the lattice structure of $H^1(\Omega)$.
 Namely, since $\Omega$ is supposed to be
a Lipschitz domain, $H^1(\Omega)$ is a {\em Dirichlet space\/}
(\cite [Theorem~5.23]{Ba-Ca}, \cite[Corollary~A.6]{Ki-St})
 in the following sense: Let
$\theta : \R \rightarrow \R$ be a uniformly Lipschitz function such that the
derivative $\theta'$ exists except at finitely many points and that
 $\theta(0) = 0$; then the induced map $\theta^*$ on $H^1(\Omega)$ 
given by $w\in H^1(\Omega)\mapsto
\theta \circ w$ is a continuous map into $H^1(\Omega)$. In particular, the map
$w \in H^1(\Omega) \mapsto w^+ = \max (0,w) = 
\frac{1}{2} (w + |w|) $ is a continuous map into $H^1(\Omega)$.

Now we are in the position to establish the following existence result.

\begin{theorem} \label{OCP-4}
Suppose that the generalized directional derivative $J^0$
satisfies the  one-sided Lipschitz condition (\ref{cc1})
and the growth condition (\ref{cc3}).
Moreover, suppose the smallness condition (\ref{cc2})
with the monotonicity constant $c_{\cal A}$ from Lemma \ref{l-coerc}.
Then there exists an optimal control 
to the optimal control problem  $(\mathcal{OCP})_4$.
\end{theorem}

{\bf Proof.\/} The proof again follows the standard pattern of existence
proofs in optimal control.

Since the cost function  $I_4$ is  bounded below, 
$$
I_4^* := \inf (\mathcal{OCP})_4 \in \R \,.
$$ 
Let $\{({\ul u}_n, {\ol u}_n)\}_{n\in\N} \subset  U_{\mbox \small {ad}}$ be a minimizing sequence of  $(\mathcal{OCP})_4$, e.g.
construct $({\ul u}_n, {\ol u}_n) \in  U_{\mbox \small {ad}}$ 
by $I_4({\ul u}_n, {\ol u}_n) < I_4^* + \frac{1}{n}$.
Since the misfit term in $I_4$ is nonnegative, 
$\|{\ul u}_n \|^2_{H^2(\Omega)}$ and $\|{\ol u}_n\|^2_{H^2(\Omega)}$ are bounded.
By reflexivity of $H^2(\Omega)$, closedness of  $U_{\mbox \small {ad}}$, and the compact embedding $H^2(\Omega)  \subset \subset H^1(\Omega)$ , 
we can pass to a subsequence of $\{ ({\ul u}_n,{\ol u}) \}_{n\in\N}$
 also denoted $\{({\ul u}_n, {\ol u}_n)\}_{n\in\N}$ such that
for some $ ({\ul u}_\infty, {\ol u}_\infty)  \in  U_{\mbox \small {ad}}$, 
\begin{eqnarray}
\label{a-4-4} \left.
\begin{array}{ll}
{\ul u}_n \rha {\ul u}_\infty \,, &  {\ol u}_n \rha {\ol u}_\infty
\qquad \mbox{ in }  H^2(\Omega) \,,\\[0.5ex]
{\ul u}_n \rightarrow {\ul u}_\infty \,, & 
 {\ol u}_n \rightarrow {\ol u}_\infty
\qquad \mbox{ in }  H^1(\Omega) \,,\\[0.5ex]
{\ul u}_n \rightarrow {\ul u}_\infty \,, & 
 {\ol u}_n \rightarrow {\ol u}_\infty
\qquad \mbox{ a.e. in }  \Omega
 \,.
\end{array}  \right\}
\end{eqnarray}
 
We claim the Mosco convergence 
${\cal C}_n \stackrel{\rm M}{\longrightarrow} {\cal C}_\infty$
for the constraint sets
\begin{eqnarray*}
&& {\cal C}_n := {\cal C}_{{\ul u}_n, {\ol u}_n} :=
\{(u,v) \in E | \, {\ul u}_n  \le u \le {\ol u}_n \mbox{ a.e. in } \Omega  \}   
\,, \\[0.5ex]
&& {\cal C}_\infty := {\cal C}_{{\ul u}_\infty, {\ol u}_\infty} :=
\{(u,v) \in E | \, {\ul u}_\infty \le u \le {\ol u}_\infty \mbox{ a.e. in } \Omega  \}  \,.
\end{eqnarray*}

To show {\bf (m1)}, let $(u_n,v_n) \in {\cal C}_n$ such that $(u_n,v_n) \rha (u,v)$ in $E$ for $n \to \infty$. Then $u_n \rha u$ in $H^1(\Omega)$.
By compact embedding $H^1(\Omega)  \subset \subset L^2(\Omega)$, 
for some subsequence $u_n \rightarrow  u$ in $ L^2(\Omega)$
and $u_n \rightarrow  u$ a.e. in $ \Omega$. Thus by (\ref{a-4-4}), 
${\ul u}_\infty \le u \le {\ol u}_\infty$ a.e. in $ \Omega$.
Hence $(u,v) \in {\cal C}$ as required and {\bf (m1)} is proven.

To show {\bf (m2)}, we exploit the above mentioned lattice structure
of $H^1(\Omega)$ and employ a cutting technique.
Let $(u,v) \in {\cal C}_\infty$.
Then ${\ul u}_\infty \le u \le {\ol u}_\infty \mbox{ a.e. in } \Omega$.
This means 
$$
\max({\ul u}_\infty , \min({\ol u}_\infty,u))
= \max({\ul u}_\infty ,u) = u \,.
$$ 
Then set 
$$u_n := \max({\ul u}_n , \min({\ol u}_n,u)), \, v_n := v \,.
$$
By construction $(u_n, v_n) \in {\cal C}_n$, moreover by (\ref{a-4-4}),
$\min({\ol u}_n,u) \rightarrow \min({\ol u}_\infty,u)$ and
$u_n \rightarrow u$ in $H^1(\Omega)$.
Hence $(u_n,v_n) \rightarrow (u,v)$ in $E$ as required. Thus {\bf (m2)} 
and the claimed Mosco convergence for the constraint sets are proved.
This entails $F_n  \stackrel{\rm M}{\longrightarrow} F_\infty$,
where $F_n := \chi_{{\cal C}_n}, F_\infty := \chi_{{\cal C}_\infty}$.
Therefore in virtue of Theorem \ref{STABLE},
${\cal S}({\ul u}_n , {\ol u}_n)  \rightarrow 
{\cal S}({\ul u}_\infty , {\ol u}_\infty )$. Thus in view of the weak 
lower semicontinuity of the norm, 
$$I_4^* \le I_4({\ul u}_\infty , {\ol u}_\infty ) 
\le \liminf_{n \to \infty} I_4({\ul u}_n , {\ol u}_n) \le I_4^* 
$$
and the theorem is proved.  
\qed 

\section{Conclusions and an Outlook} 

This paper has shown how various techniques from different fields of
mathematical analyis  can be combined to arrive at well-posedness results
for a nonlinear interface problem that models nonmonotone frictional
contact of elastic infinite media. 
In particular, we established a stability result for extended real-valued
hemivariational inequalities that extends and considerably improves the stability result of \cite{Sof2018}. Furthermore 
we investigated various  optimal control problems governed by  the interface problem and governed by a related obstacle interface problem. Based on our stability result we could present an unified
approach to existence results for optimal controls in these control problems.

In this paper we dealt with the primal HVI formulation of the
underlying interface problem in the optimal control problem. Here mixed variational formulations, see e.g. \cite{Gwi-17,SoMa2019,LiMi2019,ChMoSa},
would be another direction of research.

%outlook:
% zhenhai Liu Mig f-g-h quasimonotone ApplMathOpt 19
% existence. Chadli, MTA to appear

The next step towards numerical treatment of such optimal control problems 
is the study of relaxation methods see e.g. \cite{PaRuRe2020},
and the derivation of optimality conditions, see e.g. \cite{PeKu2018}.

A major challenge is to arrive at efficient and reliable numerical
solution methods, known in optimal control with LINEAR elliptic boundary
value problems, see e.g. \cite{ApStWi2016}.

\bibliographystyle{amsplain}

\bibliography{biblio-OC-HVIs-Interface-via-BIE}

\providecommand{\bysame}{\leavevmode\hbox to3em{\hrulefill}\thinspace}
\providecommand{\MR}{\relax\ifhmode\unskip\space\fi MR }
\providecommand{\MRhref}[2]{%
  \href{http://www.ams.org/mathscinet-getitem?mr=#1}{#2}
}
\providecommand{\href}[2]{#2}
\begin{thebibliography}{10}

\bibitem{ALY1998}
D.~R. Adams, S.~M. Lenhart, and J.~Yong, \emph{Optimal control of the obstacle
  for an elliptic variational inequality}, Appl. Math. Optim. \textbf{38}
  (1998), no.~2, 121--140.

\bibitem{ApStWi2016}
Thomas Apel, Olaf Steinbach, and Max Winkler, \emph{Error estimates for
  {N}eumann boundary control problems with energy regularization}, J. Numer.
  Math. \textbf{24} (2016), no.~4, 207--233.

\bibitem{AttBut-2014}
H.~Attouch, G.~Buttazzo, and G.~Michaille, \emph{Variational analysis in
  {S}obolev and {BV} spaces}, second ed., MOS-SIAM Series on Optimization,
  vol.~17, Society for Industrial and Applied Mathematics (SIAM), Philadelphia,
  PA; Mathematical Optimization Society, Philadelphia, PA, 2014.

\bibitem{Ba-Ca}
C.~Baiocchi and A.~Capelo, \emph{Variational and quasivariational inequalities
  - applications to free boundary problems}, John Wiley \& Sons, Inc., New
  York, 1984.

\bibitem{Ba1984}
V.~Barbu, \emph{Optimal control of variational inequalities}, Research Notes in
  Mathematics, vol. 100, Pitman (Advanced Publishing Program), Boston, MA,
  1984.

\bibitem{BlOe94}
E.~Blum and W.~Oettli, \emph{From optimization and variational inequalities to
  equilibrium problems}, Math. Student \textbf{63} (1994), no.~1-4, 123--145.

\bibitem{Capa-2014}
A.~Capatina, \emph{Variational {I}nequalities and {F}rictional {C}ontact
  {P}roblems}, Springer, Cham, 2014.

\bibitem{CCJG-1997}
C.~Carstensen and J.~Gwinner, \emph{F{EM} and {BEM} coupling for a nonlinear
  transmission problem with {S}ignorini contact}, SIAM J. Numer. Anal.
  \textbf{34} (1997), no.~5, 1845--1864.

\bibitem{CHNZ2022}
Jinxia Cen, Tahar Haddad, Van~Thien Nguyen, and Shengda Zeng,
  \emph{Simultaneous distributed-boundary optimal control problems driven by
  nonlinear complementarity systems}, J. Global Optim. \textbf{84} (2022),
  783--805.

\bibitem{ChMoSa}
O.~Chadli, R.N. Mohapatra, and B.K. Sahu, \emph{Br\'{e}zis-pseudomonotone mixed
  equilibrium problems involving a set-valued mapping with application},
  Minimax Theory Appl. (2022), to appear.

\bibitem{Clarke}
F.H. Clarke, \emph{Optimization and nonsmooth analysis}, second ed., Classics
  in Applied Mathematics, vol.~5, Society for Industrial and Applied
  Mathematics (SIAM), Philadelphia, PA, 1990.

\bibitem{DLR2011}
Juan~Carlos De~Los~Reyes, \emph{Optimal control of a class of variational
  inequalities of the second kind}, SIAM J. Control Optim. \textbf{49} (2011),
  no.~4, 1629--1658.

\bibitem{Fr1986}
A.~Friedman, \emph{Optimal control for variational inequalities}, SIAM J.
  Control Optim. \textbf{24} (1986), no.~3, 439--451.

\bibitem{GhIo2009}
M.~Ghergu and I.R. Ionescu, \emph{Structure–soil–structure coupling in
  seismic excitation and “city effect”}, Internat. J. Engrg. Sci.
  \textbf{47} (2009), no.~3, 342--354.

\bibitem{GMS-2021}
H.~Gimperlein, M.~Maischak, and E.P. Stephan, \emph{{FE}-{BE} coupling for a
  transmission problem involving microstructure}, J. Appl. Numer.Optim.
  \textbf{3} (2021), no.~2, 315--331.

\bibitem{GoeMot-2003}
D.~Goeleven and D.~Motreanu, \emph{Variational and hemivariational
  inequalities: theory, methods and applications. {V}ol. {II}}, Nonconvex
  Optimization and its Applications, vol.~70, Kluwer Academic Publishers,
  Boston, MA, 2003, Unilateral problems.

\bibitem{GMDR-2003}
D.~Goeleven, D.~Motreanu, Y.~Dumont, and M.~Rochdi, \emph{Variational and
  hemivariational inequalities: theory, methods and applications. {V}ol. {I}},
  Nonconvex Optimization and its Applications, vol.~69, Kluwer Academic
  Publishers, Boston, MA, 2003, Unilateral analysis and unilateral mechanics.

\bibitem{Gwi-81}
J.~Gwinner, \emph{On fixed points and variational inequalities---a circular
  tour}, Nonlinear Anal. \textbf{5} (1981), no.~5, 565--583.

\bibitem{Gwi-17}
\bysame, \emph{Lagrange multipliers and mixed formulations for some inequality
  constrained variational inequalities and some nonsmooth unilateral problems},
  Optimization \textbf{66} (2017), no.~8, 1323--1336.

\bibitem{GJKS18}
J.~Gwinner, B.~Jadamba, A.~A. Khan, and M.~Sama, \emph{Identification in
  variational and quasi-variational inequalities}, J. Convex Anal. \textbf{25}
  (2018), 545--569.

\bibitem{GwiSte-2018}
J.~Gwinner and E.P. Stephan, \emph{Advanced boundary element methods}, Springer
  Series in Computational Mathematics, vol.~52, Springer, Cham, 2018, Treatment
  of boundary value, transmission and contact problems.

\bibitem{Gwi-95}
Joachim Gwinner, \emph{Stability of monotone variational inequalities with
  various applications}, Variational inequalities and network equilibrium
  problems ({E}rice, 1994), Plenum, New York, 1995, pp.~123--142.

\bibitem{Gwi-22}
\bysame, \emph{From the {F}an-{KKM} principle to extended real-valued
  equilibria and to variational-hemivariational inequalities with application
  to nonmonotone contact problems}, Fixed Point Theory Algorithms Sci. Eng.
  (2022), Paper No. 4, 28.

\bibitem{HaPa1995}
J.~Haslinger and P.~D. Panagiotopoulos, \emph{Optimal control of systems
  governed by hemivariational inequalities. {E}xistence and approximation
  results}, Nonlinear Anal. \textbf{24} (1995), no.~1, 105--119.

\bibitem{HsWe-2008}
G.C. Hsiao and W.L.Wendland, \emph{Boundary integral equations}, Applied
  Mathematical Sciences, vol. 164, Springer, Berlin, 2008.

\bibitem{ItKu2000}
K.~Ito and K.~Kunisch, \emph{Optimal control of elliptic variational
  inequalities}, Appl. Math. Optim. \textbf{41} (2000), no.~3, 343--364.

\bibitem{Ki-St}
D.~Kinderlehrer and G.~Stampacchia, \emph{An introduction to variational
  inequalities and their applications}, Classics in Applied Mathematics,
  vol.~31, Society for Industrial and Applied Mathematics (SIAM), Philadelphia,
  PA, 2000.

\bibitem{Li1971}
J.-L. Lions, \emph{Optimal control of systems governed by partial differential
  equations}, Die Grundlehren der mathematischen Wissenschaften, Band 170,
  Springer-Verlag, New York-Berlin, 1971, Translated from the French by S. K.
  Mitter.

\bibitem{LYZe2021}
Jinjie Liu, Xinmin Yang, and Shengda Zeng, \emph{Optimal control and
  approximation for elliptic bilateral obstacle problems}, Commun. Nonlinear
  Sci. Numer. Simul. \textbf{102} (2021), Paper No. 105938, 17.

\bibitem{LiMi2019}
Z.~Liu, St. Mig\'{o}rski, and B.~Zeng, \emph{Existence results and optimal
  control for a class of quasi mixed equilibrium problems involving the
  {$(f,g,h)$}-quasimonotonicity}, Appl. Math. Optim. \textbf{79} (2019), no.~2,
  257--277.

\bibitem{MaiSte-2005}
M.~Maischak and E.P. Stephan, \emph{A {FEM}-{BEM} coupling method for a
  nonlinear transmission problem modelling {C}oulomb friction contact}, Comput.
  Methods Appl. Mech. Engrg. \textbf{194} (2005), no.~2-5, 453--466.

\bibitem{MiPu1984}
F.~Mignot and J.-P. Puel, \emph{Optimal control in some variational
  inequalities}, SIAM J. Control Optim. \textbf{22} (1984), no.~3, 466--476.

\bibitem{Ochal}
S.~Mig\'{o}rski, A.~Ochal, and M.~Sofonea, \emph{Nonlinear inclusions and
  hemivariational inequalities}, Advances in Mechanics and Mathematics,
  vol.~26, Springer, New York, 2013, Models and analysis of contact problems.

\bibitem{Mos1969}
U.~Mosco, \emph{Convergence of convex sets and of solutions of variational
  inequalities}, Advances in Math. \textbf{3} (1969), 510--585.

\bibitem{Naniewicz}
Z.~Naniewicz and P.D. Panagiotopoulos, \emph{Mathematical theory of
  hemivariational inequalities and applications}, Monographs and Textbooks in
  Pure and Applied Mathematics, vol. 188, Marcel Dekker, Inc., New York, 1995.

\bibitem{OvGw-Rassias-14}
Gwinner~J. Ovcharova, N., \emph{On the discretization of pseudomonotone
  variational inequalities with an application to the numerical solution of the
  nonmonotone delamination problem}, Optimization in Science and Engineering,
  Springer, New York, 2014, pp.~393--405.

\bibitem{Ov2017}
N.~Ovcharova, \emph{On the coupling of regularization techniques and the
  boundary element method for a hemivariational inequality modelling a
  delamination problem}, Math. Methods Appl. Sci. \textbf{40} (2017), no.~1,
  60--77.

\bibitem{OvGw-2014}
N.~Ovcharova and J.~Gwinner, \emph{A study of regularization techniques of
  nondifferentiable optimization in view of application to hemivariational
  inequalities}, J. Optim. Theory Appl. \textbf{162} (2014), no.~3, 754--778.

\bibitem{Pan-1993}
P.D. Panagiotopoulos, \emph{Hemivariational inequalities}, Springer-Verlag,
  Berlin, 1993, Applications in mechanics and engineering.

\bibitem{PaRuRe2020}
Nikolaos~S. Papageorgiou, Vicen\c{t}iu~D. R\u{a}dulescu, and Du\v{s}an~D.
  Repov\v{s}, \emph{Relaxation methods for optimal control problems}, Bull.
  Math. Sci. \textbf{10} (2020), no.~1, 2050004, 24.

\bibitem{PaWi2018}
Nikolaos~S. Papageorgiou and Patrick Winkert, \emph{Applied nonlinear
  functional analysis}, De Gruyter Graduate, De Gruyter, Berlin, 2018, An
  introduction.

\bibitem{Pa1977}
F.~Patrone, \emph{On the optimal control for variational inequalities}, J.
  Optim. Theory Appl. \textbf{22} (1977), no.~3, 373--388.

\bibitem{PeKu2018}
Zijia Peng and Karl Kunisch, \emph{Optimal control of elliptic
  variational-hemivariational inequalities}, J. Optim. Theory Appl.
  \textbf{178} (2018), no.~1, 1--25.

\bibitem{RTR}
R.~Tyrrell Rockafellar, \emph{Convex analysis}, Princeton Landmarks in
  Mathematics, Princeton University Press, Princeton, NJ, 1997, Reprint of the
  1970 original, Princeton Paperbacks.

\bibitem{SaSch2011}
S.A. Sauter and Ch. Schwab, \emph{Boundary element methods}, Springer Series in
  Computational Mathematics, vol.~39, Springer-Verlag, Berlin, 2011.

\bibitem{scholz2019}
C.H. Scholz, \emph{The mechanics of earthquakes and faulting}, Cambridge
  {U}niversity {P}ress, 2019.

\bibitem{SofMig-2018}
M.~Sofonea and S.~Mig\'{o}rski, \emph{Variational-hemivariational inequalities
  with applications}, Monographs and Research Notes in Mathematics, CRC Press,
  Boca Raton, FL, 2018.

\bibitem{Sof2018}
Mircea Sofonea, \emph{Convergence results and optimal control for a class of
  hemivariational inequalities}, SIAM J. Math. Anal. \textbf{50} (2018), no.~4,
  4066--4086.

\bibitem{Sof2019}
\bysame, \emph{Optimal control of a class of variational-hemivariational
  inequalities in reflexive {B}anach spaces}, Appl. Math. Optim. \textbf{79}
  (2019), no.~3, 621--646.

\bibitem{SoMa2019}
Mircea Sofonea, Andaluzia Matei, and Yi-bin Xiao, \emph{Optimal control for a
  class of mixed variational problems}, Z. Angew. Math. Phys. \textbf{70}
  (2019), no.~4, Paper No. 127, 17.

\bibitem{Tr2010}
Fredi Tr\"{o}ltzsch, \emph{Optimal control of partial differential equations},
  Graduate Studies in Mathematics, vol. 112, American Mathematical Society,
  Providence, RI, 2010, Theory, methods and applications, Translated from the
  2005 German original by J\"{u}rgen Sprekels.

\bibitem{wang2018soil}
G.~Wang, M.~Yuan, Y.~Miao, J.~Wu, and Y.~Wang, \emph{Experimental study on
  seismic response of underground tunnel-soil-surface structure interaction
  system}, Tunnelling and Underground Space Technology \textbf{76} (2018),
  145--159.

\bibitem{Zei-2}
E.~Zeidler, \emph{Nonlinear functional analysis and its applications.
  {II}/{B}}, Springer-Verlag, New York, 1990, Nonlinear monotone operators,
  Translated from the German by the author and Leo F. Boron.

\end{thebibliography}

\end{document}